\documentclass[12pt,a4paper]{article}
\usepackage{amsmath,amssymb,amsthm,mathtools}
\usepackage{mathrsfs}
\usepackage{setspace,graphicx,array,tikz,url}
\usepackage{cite}
\usepackage[hidelinks]{hyperref}
\usepackage{fancyhdr}
\usepackage{lastpage}
\usepackage[top=3cm,left=3cm,right=3cm,bottom=3cm]{geometry}
\usepackage{enumitem}
\usepackage{titlesec}
\usepackage{etoolbox}
\usepackage{microtype}

\patchcmd{\thebibliography}
  {\advance\leftmargin\labelsep}
  {\labelsep=.5em\leftmargin=2.5em\itemindent=0em\advance\leftmargin\labelsep}
  {}{}

\titlespacing*{\section}{0pt}{3em}{2em}

\newtheoremstyle{spaced}
  {15pt}   
  {20pt}   
  {\itshape} 
  {}       
  {\bfseries} 
  {.}      
  {.5em}   
  {}       

\theoremstyle{spaced}

\newtheorem{theorem}{Theorem}[section]
\newtheorem{definition}{Definition}[section]

\newtheorem{lemma}{Lemma}[section]

\newtheorem{corollary}{Corollary}[section]

\title{Isomorphism in Union-Closed Sets}
\author{
  Mohammad Javad Moghaddas Mehr\thanks{m.moghadas11235@gmail.com}}

\date{January 5, 2025}

\begin{document}
\maketitle
\vspace{5em}

\begin{abstract}
	We prove that for any isomorphism \(h: \mathcal{K}_1 \to \mathcal{K}_2\) between pure union-closed families, there exists a hyperisomorphism \(H: \bigcup \mathcal{K}_1 \to \bigcup \mathcal{K}_2\) such that \(h(A) = \{ H(a) \mid a \in A \}\), for all \(A \in \mathcal{K}_1\). Since every union-closed family forms a lattice under inclusion, this result establishes a strong connection between the two frameworks. More precisely, any such family can be uniquely reconstructed from its lattice up to isomorphism. Hence, the lattice representation provides a faithful encoding, offering a perspective that may yield new insights into problems on union-closed families, including Frankl’s union-closed sets conjecture.
\end{abstract}


\section{Introduction}
\paragraph{}
The motivation for this paper arose while studying structural properties of
union-closed families, which revealed a strong connection between union-closed families
and lattices.
It is straightforward to show that any union-closed family containing the empty set can be uniquely represented as a lattice under set inclusion. Together with
Theorem~\ref{thm:main}, this suggests that a lattice-theoretic approach
can be highly effective in investigating properties of union-closed families. More
precisely, converting a union-closed family into a lattice structure preserves all
essential information about the family. \\

One of the central open questions concerning union-closed families is Frankl’s
conjecture~\cite{frankl1995}, which states that for any finite union-closed family
there exists an element that belongs to at least half of the sets in the family.
Despite decades of effort, the conjecture remains unresolved; see, for example,~\cite{Bruhn2015,moghaddasmehr2023,gilmer2022,DasWu2024,Cambie2022,Yu2023,Alweiss2024,Karpas2017,Colbert2024,Gendler2025,LuRaz2024,CarvalhoMachiavelo2024} for background and recent progress.\\

In fact, the lattice formulation of Frankl’s conjecture has already been established
for several important classes. Rival first observed that the conjecture holds for
distributive and geometric lattices (without proof)~\cite{Rival1985}, and Poonen later
gave a complete proof and extended the result to relatively complemented
lattices~\cite{Poonen1992}. Abe and Nakano proved the conjecture for modular
lattices~\cite{AbeNakano1998}, while Reinhold established it for lower semimodular
lattices~\cite{Reinhold2000}, which remains the strongest known result for a standard
lattice class. Abe and Nakano also showed it for lower quasi-semimodular
lattices~\cite{AbeNakano2000}, and Czédli and Schmidt confirmed it for large
semimodular lattices and for planar semimodular lattices~\cite{CzedliSchmidt2008}.
Joshi, Waphare, and Kavishwar proved the conjecture for dismantlable lattices, from which they derived as corollaries the cases of upper semimodular lattices of breadth at most two~\cite{JoshiWaphareKavishwar2016}. Later, Joshi and Waphare explicitly
established the conjecture for all lattices of breadth two~\cite{JoshiWaphare2019}.
Independently, Abdollahi, Woodroofe, and Zaimi proved the conjecture for subgroup
lattices of finite groups and, more generally, for lattices with a modular coatom—a
class that includes supersolvable and dually semimodular
lattices~\cite{AbdollahiWoodroofeZaimi2017}. \\

In the next section, we introduce the basic concepts needed for the proof of Theorem~\ref{thm:main}, which states that for every
\textit{\hyperref[def:isomorphism]{isomorphism}} \(h: \mathcal{K}_1 \to \mathcal{K}_2\) between two
\textit{\hyperref[def:pure_collection]{pure union-closed families of sets}}, there exists a corresponding
\textit{\hyperref[def:hyperisomorphism]{hyperisomorphism}}
\(H: \bigcup \mathcal{K}_1 \to \bigcup \mathcal{K}_2\) such that
\(
h(A) = \{ H(a) \mid a \in A \}, \quad \forall\, A \in \mathcal{K}_1.
\)

\section{Preliminaries}
\paragraph{}
For an integer \(n \in \mathbb{N}\), we define \([n] := \{k \in \mathbb{N} \mid k \leq n\}\). A family \(\mathcal{K} \subseteq 2^{[n]}\) is called \textbf{union-closed} if for all \(A, B \in \mathcal{K}\), it holds that \(A \cup B \in \mathcal{K}\). A \textbf{lattice} is a poset \((L, \leq)\) in which any two elements of \(L\) have a unique \textit{meet} and a unique \textit{join}.

\begin{definition}[Homomorphism]\label{def:homomorphism}
	Let \( \mathcal{K}_1 \) and \( \mathcal{K}_2 \) be union-closed families of sets.
	A mapping \( h: \mathcal{K}_1 \to \mathcal{K}_2 \) is called a \textbf{homomorphism} if, for all
	\( A_1, A_2 \in \mathcal{K}_1 \), the following holds:
	\[
		h(A_1 \cup A_2) = h(A_1) \cup h(A_2).
	\]
\end{definition}

\begin{definition}[Isomorphism]\label{def:isomorphism}
	Let \( \mathcal{K}_1 \) and \( \mathcal{K}_2 \) be union-closed families of sets.
	A homomorphism \( h: \mathcal{K}_1 \to \mathcal{K}_2 \) is called an \textbf{isomorphism} if it is bijective.
\end{definition}

\paragraph{}
It is straightforward to see that any union-closed family containing the empty set forms a lattice under inclusion. In this way, homomorphisms and isomorphisms of families correspond exactly to homomorphisms and isomorphisms of the associated lattices.

\begin{lemma}
	Let \( \mathcal{K}_1 \) and \( \mathcal{K}_2 \) be two union-closed families, and \( h: \mathcal{K}_1 \to \mathcal{K}_2 \) be a \textit{homomorphism} between them. If \( A_1, A_2 \in \mathcal{K}_1 \) such that \( A_1 \subseteq A_2 \), then:
	\[
		h(A_1) \subseteq h(A_2).
	\]
\end{lemma}

\begin{proof}
	By the homomorphism property, we have:
	\[
		h(A_2) = h(A_1 \cup A_2) = h(A_1) \cup h(A_2).
	\]
	Thus, \( h(A_1) \subseteq h(A_2) \).
\end{proof}

\begin{corollary}\label{cor:subset}
	Let \( \mathcal{K}_1 \) and \( \mathcal{K}_2 \) be two union-closed families, and \( h: \mathcal{K}_1 \to \mathcal{K}_2 \) be an \textit{isomorphism} between them.  If \( A_1, A_2 \in \mathcal{K}_1 \) such that \( A_1 \subset A_2 \), then:
	\[
		h(A_1) \subset h(A_2).
	\]
\end{corollary}

\begin{lemma}
	Let \(\mathcal{K}\) be a union-closed family of sets, and \( h: \mathcal{K} \to 2^{[n]} \) be a homomorphism. Then the image of \( \mathcal{K} \) under \( h \) forms a union-closed family of sets.
\end{lemma}

\begin{proof}
	Let \( B_1, B_2 \in h(\mathcal{K}) \). Then there exist \( A_1, A_2 \in \mathcal{K} \) such that \( B_1 = h(A_1) \) and \( B_2 = h(A_2) \). Therefore:
	\[
		B_1 \cup B_2 = h(A_1) \cup h(A_2) = h(A_1 \cup A_2).
	\]
	Thus, \( B_1 \cup B_2 \in h(\mathcal{K}) \), proving that \( h(\mathcal{K}) \) is union-closed.
\end{proof}

\begin{corollary}
	Let \(\mathcal{K}\) be a union-closed family of sets, and \( h: \mathcal{K} \to 2^{[n]} \) be
	an injective homomorphism. Then \( h: \mathcal{K} \to h(\mathcal{K}) \) is an isomorphism.
\end{corollary}

\begin{definition}[Redundant Element]
	Let \( \mathcal{K} \subseteq 2^{[n]} \). An element \( z \in \bigcup \mathcal{K} \) is called \textbf{redundant}, if removing \( z \) from every set in \( \mathcal{K} \) does not reduce the cardinality of the collection. Specifically, we define:
	\[
		\mathcal{K}^{\setminus z} = \{X \setminus \{z\} \mid X \in \mathcal{K}\},
	\]
	where \( |\mathcal{K}| = |\mathcal{K}^{\setminus z}| \). The collection \( \mathcal{K}^{\setminus z} \) is called the \textbf{reduced collection}.
\end{definition}

\begin{definition}[Pure Collection]\label{def:pure_collection}
	A collection \(\mathcal{K} \subseteq 2^{[n]}\) is called \textbf{pure} if it does not have any \textbf{redundant} element.
\end{definition}

\paragraph{}
The next step is to remove redundant elements from a collection. Eliminating one element can change whether others remain redundant, so the process must be carried out iteratively until no redundant elements remain. The result of this process may depend on the order in which elements are removed. Corollary~\ref{cor:pure} shows, however, that all such outcomes are isomorphic. We therefore speak of the \textbf{purified collection} of \(\mathcal{K}\), denoted by \(\mathcal{K}^*\). If \(\mathcal{K}\) is already pure, then clearly \(\mathcal{K} = \mathcal{K}^*\).

\begin{lemma}
	Let \( \mathcal{K} \) be a union-closed family of sets, and \( z \) be a redundant element of \( \mathcal{K} \). Then there exists an isomorphism between \( \mathcal{K} \) and \( \mathcal{K}^{\setminus z} \).
\end{lemma}

\begin{proof}
	Define \( h: \mathcal{K} \to \mathcal{K}^{\setminus z} \) by \( h(A) = A \setminus \{z\} \). Let \( A, B \in \mathcal{K} \). Then:
	\[
		h(A \cup B) = (A \cup B) \setminus \{z\} = (A \setminus \{z\}) \cup (B \setminus \{z\}) = h(A) \cup h(B).
	\]
	This shows that \( h \) preserves the union operation, so \( h \) is a homomorphism.
	It is clear that \( h \) is surjective because every element of \( \mathcal{K}^{\setminus z} \) is of the form \( A \setminus \{z\} \) for some \( A \in \mathcal{K} \).
	To show that \( h \) is injective, suppose \( h(A) = h(B) \) for some \( A, B \in \mathcal{K} \). This implies:
	\[
		A \setminus \{z\} = B \setminus \{z\}.
	\]

	If \( A \neq B \), removing \( z \) would change the size of \( \mathcal{K} \), contradicting the assumption that \( z \) is a redundant element. Thus, \( A = B \), which implies that \( h \) is injective. Since \( h \) is both a homomorphism and a bijection, it follows that \( h \) is an isomorphism between \( \mathcal{K} \) and \( \mathcal{K}^{\setminus z} \).
\end{proof}

\begin{corollary}\label{cor:pure}
	For any union-closed family \( \mathcal{K} \), there exists an isomorphism \( h: \mathcal{K} \to \mathcal{K}^* \).
\end{corollary}

\section{Cardinality Theorem}

\paragraph{}
This section establishes the \textit{\hyperref[thm:cardinality]{Cardinality Theorem}}, which we will use in the proof of the main theorem. We first record some terminology and supporting lemmas.

\begin{definition}
	Let \(\mathcal{K} \subseteq 2^{[n]}\) be a poset under set inclusion. An element \( X \in \mathcal{K} \) is called \textbf{minimal} if it has no proper subset in \(\mathcal{K}\). The set of all minimal elements in \(\mathcal{K}\) is denoted by \(\mathcal{K}^{\bot}\) and is defined as:
	\[
		\mathcal{K}^{\bot} = \{ X \in \mathcal{K} \mid \nexists A \in \mathcal{K} \text{ such that } A \subset X \}.
	\]
\end{definition}

\begin{figure}
	\centering
	\begin{minipage}{0.45\textwidth}
		\centering
		\begin{tikzpicture}[scale=1, transform shape, every node/.style={draw, circle, minimum size=0.5cm}]
			\node[fill=red, opacity=0.5, draw=red] (empty) at (0, 0) {};
			\node (x) at (-2, 2) {};
			\node (z) at (2, 2) {};
			\node (xy) at (-2, 4) {};
			\node (xz) at (0, 4) {};
			\node (yz) at (2, 4) {};
			\node (xyz) at (0, 6) {};

			\draw (empty) -- (x);
			\draw (empty) -- (z);
			\draw (x) -- (xy);
			\draw (x) -- (xz);
			\draw (z) -- (xz);
			\draw (z) -- (yz);
			\draw (xy) -- (xyz);
			\draw (xz) -- (xyz);
			\draw (yz) -- (xyz);
		\end{tikzpicture}
	\end{minipage}
	\hfill
	\begin{minipage}{0.45\textwidth}
		\centering
		\begin{tikzpicture}[scale=1, transform shape, every node/.style={draw, circle, minimum size=0.5cm}]
			\node[fill=red, opacity=0.5, draw=red] (y) at (2, 0) {};
			\node[fill=red, opacity=0.5, draw=red] (x) at (-2, 2) {};
			\node (z) at (2, 2) {};
			\node (xy) at (-2, 4) {};
			\node (xz) at (0, 4) {};
			\node (yz) at (2, 4) {};
			\node (xyz) at (0, 6) {};

			\draw (y) -- (z);
			\draw (x) -- (xy);
			\draw (x) -- (xz);
			\draw (z) -- (xz);
			\draw (z) -- (yz);
			\draw (xy) -- (xyz);
			\draw (xz) -- (xyz);
			\draw (yz) -- (xyz);
		\end{tikzpicture}
	\end{minipage}
	\caption{Two lattice diagrams; the red nodes are the elements of \(\mathcal{K}^{\bot}\).}\label{fig:lattice_with_rank}
\end{figure}

\begin{lemma}\label{lem:uc}
	Let \( \mathcal{K} \) be a union-closed family, and \( X \in \mathcal{K}^{\bot} \). Then \( \mathcal{K} \setminus \{X\} \) is also union-closed.
\end{lemma}

\begin{proof}
	Since \(X\) is minimal, no member of \(\mathcal{K}\setminus\{X\}\) is contained in \(X\); hence if \(A,B\in\mathcal{K}\setminus\{X\}\) then \(A\cup B\neq X\). Because \(\mathcal{K}\) is union-closed, \(A\cup B\in\mathcal{K}\), and therefore \(A\cup B\in\mathcal{K}\setminus\{X\}\).
\end{proof}

\newpage

\begin{lemma}\label{lem:minimality}
	Let \( \mathcal{K}_1 \) and \( \mathcal{K}_2 \) be union-closed families, and \( h: \mathcal{K}_1 \to \mathcal{K}_2 \) be an isomorphism between them. If \( X \) is a minimal element of \( \mathcal{K}_1 \), then \( h(X) \) is a minimal element of \( \mathcal{K}_2 \).
\end{lemma}

\begin{proof}
	Assume, for contradiction, that \( h(X) \) is not minimal in \( \mathcal{K}_2 \). Then there exists \( B \in \mathcal{K}_2 \) with \( B \subset h(X) \). Since \( h \) is surjective, there is some \( A \in \mathcal{K}_1 \) such that \( h(A) = B \). Clearly \( A \neq X \).
	Now consider
	\[
		h(A \cup X) = h(A) \cup h(X) = B \cup h(X) = h(X).
	\]
	Because \( h \) is injective, it follows that \( A \cup X = X \), hence \( A \subset X \). This contradicts the minimality of \( X \) in \( \mathcal{K}_1 \). Therefore, \( h(X) \) must be minimal in \( \mathcal{K}_2 \).
\end{proof}

\begin{lemma}\label{lem:chain}
	Let \(\mathcal{K}\) be a pure union-closed family of sets. If \(|\mathcal{K}^{\bot}| = 1\), then \(\mathcal{K}^{\bot} = \{\varnothing\}\).
\end{lemma}

\begin{proof}
	By assumption, there exists exactly one set, denoted by \( A \), that belongs to \(\mathcal{K}^{\bot}\). It is straightforward to verify that every maximal chain in \(\mathcal{K}\) must start from a member of \(\mathcal{K}^{\bot}\). In this case, every maximal chain can be represented as:
	\[
		A \subset C_1 \subset C_2 \subset \cdots \subset \bigcup \mathcal{K}.
	\]
	Since every element \( X \in \mathcal{K} \) lies on a maximal chain, it follows that \( A \subset X \). If \( A \) contained any elements, those elements would be redundant, as \( A \) is a subset of all elements in \(\mathcal{K}\). However, because \(\mathcal{K}\) is a pure collection of sets, \( A \) must be the empty set.
\end{proof}

\begin{theorem}[Cardinality]\label{thm:cardinality}
	Let \( \mathcal{K}_1 \) and \( \mathcal{K}_2 \) be two \textit{\hyperref[def:pure_collection]{pure union-closed families of sets}}. If there exists an isomorphism \( h: \mathcal{K}_1 \to \mathcal{K}_2 \), then for every \( A \in \mathcal{K}_1 \), we have:
	\[
		|A| = |h(A)|.
	\]
\end{theorem}

\newpage

\begin{proof}
	We prove the statement by induction.

	\paragraph{}
	\textbf{Base Case:} Consider \( |\mathcal{K}_1| = |\mathcal{K}_2| = 1 \). In this case, \( \mathcal{K}_1 = \{A\} \) and \( \mathcal{K}_2 = \{B\} \) for some \( A, B \in 2^{[n]} \). Both \( A \) and \( B \) must necessarily be empty sets. If \( A \) or \( B \) were not empty, their elements would be redundant, and removing those elements would not affect the cardinality of \( \mathcal{K}_1 \) or \( \mathcal{K}_2 \).	Thus, for the base case, we have \( |\bigcup \mathcal{K}_1| = |\bigcup \mathcal{K}_2| = 0 \), as desired.

	\paragraph{}
	\textbf{Inductive Step:} Assume that the theorem holds for all pure union-closed families with
	cardinality strictly less than \( n \). Now, consider \( \mathcal{K}_1 \) and
	\( \mathcal{K}_2 \) such that \( |\mathcal{K}_1| = |\mathcal{K}_2| = n \), and let
	\( h: \mathcal{K}_1 \to \mathcal{K}_2 \) be an isomorphism.
	Let \( X \in \mathcal{K}_1^{\bot} \). By Lemma~\ref{lem:minimality}, we have
	\( h(X) = Y \) where \( Y \in \mathcal{K}_2^{\bot} \). Consider the restriction
	\( h: \mathcal{K}_1 \setminus \{X\} \to \mathcal{K}_2 \setminus \{Y\} \).
	By Lemma~\ref{lem:uc}, this restricted mapping remains an isomorphism.

	\paragraph{Case 1:}
	Suppose \( \mathcal{K}_1 \setminus \{X\} \) and \( \mathcal{K}_2 \setminus \{Y\} \) are pure collections. By the induction hypothesis, for all \( A \in \mathcal{K}_1 \setminus \{X\} \), we have \( |A| = |h(A)| \). If \( |\mathcal{K}_1^{\bot}| = 1 \), then by Lemma~\ref{lem:chain}, \( X = \varnothing \) and \( Y = \varnothing \) as well. Consequently, for all \( A \in \mathcal{K}_1 \), we conclude \( |A| = |h(A)| \).

	\paragraph{}
	On the other hand, if \( |\mathcal{K}_1^{\bot}| > 1 \), let \( Z \in \mathcal{K}_1^{\bot} \) such that \( Z \neq X \). If \( \mathcal{K}_1 \setminus \{Z\} \) and \( \mathcal{K}_2 \setminus \{h(Z)\} \) are pure collections, then by the induction hypothesis, we again obtain \( |X| = |Y| \), completing the proof for this case.

	\paragraph{}
	Thus, the theorem is proved for the scenario where removing one element of \( \mathcal{K}_1^{\bot} \) and its corresponding map results in pure collections.
	Next, we consider the case where there exists at least one element in \( \mathcal{K}_1^{\bot} \) or \( \mathcal{K}_2^{\bot} \) such that its removal results in a non-pure collection.

	\paragraph{Case 2:}
	Without loss of generality, suppose \( \mathcal{K}_1 \setminus \{X\} \) is not pure and contains a redundant element, denoted by \( z \). To proceed, we establish three key facts:
	\begin{enumerate}
		\item \( z \in C \) for all \( C \in \mathcal{K}_1 \setminus \{X\} \).
		\item \( \mathcal{K}_1 \setminus \{X\} \) has at most one redundant element.
		\item If \( \mathcal{K}_1 \setminus \{X\} \) has a redundant element, then \( \mathcal{K}_2 \setminus \{Y\} \) also has a redundant element.
	\end{enumerate}

	\paragraph{}
	First, consider \( z \) as a redundant element of \( \mathcal{K}_1 \setminus \{X\} \) but not of \( \mathcal{K}_1 \). This implies there exists exactly one element \( A_X \in \mathcal{K}_1 \setminus \{X\} \) such that \( A_X \setminus \{z\} = X \), and no distinct \( A_1, A_2 \in \mathcal{K}_1 \setminus \{X\} \) satisfy \( A_1 \setminus \{z\} = A_2 \). Now, for the sake of contradiction, suppose there exists \( B \in \mathcal{K}_1 \setminus \{X\} \) such that \( z \notin B \). Then \( z \notin B \cup X \). Let \( A_1 = A_X \cup B \) and \( A_2 = B \cup X \). It is straightforward to verify that \( A_1 \setminus \{z\} = A_2 \), leading to a contradiction. Thus, if \( z \) is a redundant element of \( \mathcal{K}_1 \setminus \{X\} \), it must belong to all its members.

	\paragraph{}
	Next, suppose \( t \) is another redundant element of \( \mathcal{K}_1 \setminus \{X\} \). From the previous fact, we know that \( t \) must belong to all members of \( \mathcal{K}_1 \setminus \{X\} \) but not to \( X \). However, \( t \) cannot satisfy this condition because \( A_X \setminus \{z\} = X \). Thus no additional redundant element can exist. Therefore, \( \mathcal{K}_1 \setminus \{X\} \) can contain at most one redundant element.

	\paragraph{}
	For the third fact, let \( z \) be the redundant element of \( \mathcal{K}_1 \setminus \{X\} \) and \( A_X = X \cup \{z\} \). Define \( R = h(A_X) \setminus Y \), which is non-empty by Corollary~\ref{cor:subset}. We show that \( R \) is a subset of all members of \( \mathcal{K}_2 \setminus \{Y\} \). For the sake of contradiction, suppose there exists \( D \in \mathcal{K}_2 \setminus \{Y\} \) such that \( R \not\subset D \). By fact 1, we know
	\(z \in h^{-1}(D)\). Therefore:
	\[
		h^{-1}(D) \cup X = h^{-1}(D) \cup A_X \implies D \cup Y = D \cup h(A_X).
	\]
	\paragraph{}
	This leads to a contradiction, since \( R \) is not a subset of \( D \cup Y \), but it is a subset of \( D \cup h(A_X) \).
	Thus, \( R \) must be a subset of all members of \( \mathcal{K}_2 \setminus \{Y\} \). Therefore, the members of \( R \) are redundant. As previously proven, \( \mathcal{K}_2 \setminus \{Y\} \) has at most one redundant element. Hence, \( |R| = 1 \), and we denote \( R = \{r\} \).

	\paragraph{}
	The mapping
	\(h^*: {(\mathcal{K}_1 \setminus \{X\})}^* \to {(\mathcal{K}_2 \setminus \{Y\})}^*\),
	defined as \(h^*(A \setminus \{z\}) = h(A) \setminus \{r\}\) where \(A \in \mathcal{K}_1 \setminus \{X\} \), is an isomorphism between two
	pure union-closed families of sets. Furthermore, the cardinalities satisfy
	\( |{(\mathcal{K}_1 \setminus \{X\})}^*| = |{(\mathcal{K}_2 \setminus \{Y\})}^*| = n-1 \).
	By the induction hypothesis, for all
	\(A \setminus \{z\} \in {(\mathcal{K}_1 \setminus \{X\})}^*\), we have:
	\[
		|A \setminus \{z\}| = |h^*(A \setminus \{z\})| = |h(A) \setminus \{r\}|.
	\]
	Since for all \(A \in \mathcal{K}_1 \setminus \{X\}\), we have \(z \in A\) and \(r \in h(A)\), it follows that \(|A| = |h(A)|\).

	\paragraph{}
	From fact 1 and fact 3, there exists a set \(A_X \in \mathcal{K}_1\) such that \(z \in A_X\) but \(z \not\in X\), and \(r \in h(A_X)\) but \(r \not\in Y\). This implies:
	\[
		|X| = |A_X| - 1 = |h(A_X)| - 1 = |Y|.
	\]
	Thus, for all \(A \in \mathcal{K}_1\), it follows that \(|A| = |h(A)|\).

\end{proof}

\paragraph{}
We conclude with two corollaries used later in the proof of the main theorem.

\begin{corollary}\label{cor:ci}
	Let \( \mathcal{K}_1 \) and \( \mathcal{K}_2 \) be two pure union-closed families of sets. If there exists an isomorphism \( h: \mathcal{K}_1 \to \mathcal{K}_2 \), then:
	\[
		\left|\bigcup \mathcal{K}_1\right| = \left|\bigcup \mathcal{K}_2\right|.
	\]
\end{corollary}

\begin{corollary}\label{cor:cj}
	Let \( \mathcal{K}_1 \) and \( \mathcal{K}_2 \) be two pure union-closed families of sets. If there exists an isomorphism \( h: \mathcal{K}_1 \to \mathcal{K}_2 \), then for any \(A, B \in \mathcal{K}_1\), the following properties hold:
	\begin{enumerate}
		\item \(|A \cup B| = |h(A) \cup h(B)|\).
		\item \(|A \cap B| = |h(A) \cap h(B)|\).
		\item \(|A \setminus B| = |h(A) \setminus h(B)|\).
		\item \(|A^c| = |h{(A)}^c|\),
	\end{enumerate}
	where \(A^c = \bigcup \mathcal{K}_1 \setminus A\) and
	\(h{(A)}^c = {\bigcup \mathcal{K}_2} \setminus h(A)\).

\end{corollary}

\section{Hyperisomorphism}
\paragraph{}
In this section we show that the internal structure of a pure union-closed family is faithfully preserved under isomorphism, in the sense that every isomorphism arises from a bijection of the ground sets.

\begin{definition}[Hyperisomorphism]\label{def:hyperisomorphism}
	Let \( \mathcal{K}_1 \) and \( \mathcal{K}_2 \) be two union-closed families of sets. A bijective mapping
	\( H: \bigcup \mathcal{K}_1 \to \bigcup \mathcal{K}_2 \) is called a \textbf{hyperisomorphism} if the induced mapping \( h: \mathcal{K}_1 \to \mathcal{K}_2 \), defined by \( h(A) = \{ H(a) \mid a \in A \} \), is an isomorphism.
\end{definition}

\paragraph{}
To establish our main theorem, we require the following lemmas. First, we introduce the notation \(\mathcal{K}^i\), defined as:
\[
	\mathcal{K}^i = \{A \in \mathcal{K} \mid i \in A\}.
\]

\begin{lemma}\label{lem:mi}
	Let \(\mathcal{K}\) be a pure union-closed family of sets. Then, for all \(i, j \in \bigcup \mathcal{K}\), the following equivalence holds:
	\[
		i = j \iff \mathcal{K}^i = \mathcal{K}^j.
	\]
\end{lemma}

\begin{proof}
	Suppose \(i = j\). By definition, it follows that \(\mathcal{K}^i = \mathcal{K}^j\).
	Conversely, assume \(\mathcal{K}^i = \mathcal{K}^j\), but \(i \neq j\). This would imply that \(i\) and \(j\) are redundant elements in \(\mathcal{K}\), contradicting the purity of \(\mathcal{K}\). Hence, it must be the case that \(i = j\).
\end{proof}

\begin{lemma}
	Let \(\mathcal{K}\) be a union-closed family of sets. Then, for any distinct elements \(i, j \in \bigcup \mathcal{K}\), the following inequality holds:
	\[
		\left| \left(\bigcap \mathcal{K}^i\right) \cup \left(\bigcap \mathcal{K}^j\right) \right| \geq 2.
	\]
\end{lemma}

\begin{proof}
	This follows since \(i\) and \(j\), being distinct, are both in the union.
\end{proof}

\begin{corollary}\label{cor:ui}
	Let \(\mathcal{K}\) be a union-closed family of sets. Then, for any distinct elements
	\(\{a_i\}_{i=1}^n \subseteq \bigcup \mathcal{K}\), the following inequality holds:
	\[
		\left| \bigcup_{i=1}^n \bigcap \mathcal{K}^{a_i} \right| \geq n.
	\]
\end{corollary}

\begin{lemma}\label{lem:l3n}
	Let \(\mathcal{K}\) be a union-closed family of sets, and \(i\) be an arbitrary element of \(\bigcup \mathcal{K}\).
	If \(\bigcap \mathcal{K}_1^i = \{a_1, \dots, a_n\}\), then:
	\[
		\left| \bigcup_{j=1}^n \bigcap \mathcal{K}_1^{a_j} \right| = n.
	\]
\end{lemma}

\begin{proof}
	Let \(a\) be an arbitrary element of \(\bigcap \mathcal{K}^i\). Since \(\mathcal{K}^i \subseteq \mathcal{K}^a\), it follows that \(\bigcap \mathcal{K}^a \subseteq \bigcap \mathcal{K}^i\). Consequently, we have:
	\[
		\bigcup_{j=1}^n \bigcap \mathcal{K}^{a_j} \subseteq \bigcap \mathcal{K}^i,
	\]
	which implies:
	\[
		\left| \bigcup_{j=1}^n \bigcap \mathcal{K}_1^{a_j} \right| \leq n.
	\]
	Applying Corollary~\ref{cor:ui}, we conclude that:
	\[
		\left| \bigcup_{j=1}^n \bigcap \mathcal{K}_1^{a_j} \right| = n.
	\]
	This completes the proof.
\end{proof}

\begin{lemma}\label{lem:mp}
	Let \(\mathcal{K}_1\) and \(\mathcal{K}_2\) be two pure union-closed families of sets, and
	\(h: \mathcal{K}_1 \to \mathcal{K}_2\) be an isomorphism between them. Then, for each \(i \in \bigcup \mathcal{K}_1\), there exists a unique \(j \in \bigcup \mathcal{K}_2\) such that:
	\[
		\mathcal{K}_2^j = h(\mathcal{K}_1^i),
	\]
	where \(h(\mathcal{K}_1^i) = \{h(A) \mid A \in \mathcal{K}_1^i \}\).
\end{lemma}

\begin{proof}
	Let \(i \in \bigcup \mathcal{K}_1\). Suppose \(\bigcap \mathcal{K}_1^i = \{a_1, \dots, a_n\}\). By Lemma~\ref{lem:l3n}, we have:
	\[
		\left| \bigcup_{j=1}^n \bigcap \mathcal{K}_1^{a_j} \right| = n.
	\]
	By Corollary~\ref{cor:cj}, we can assume \(\bigcap h(\mathcal{K}_1^i) = \{b_1, \dots, b_n\}\).

	\paragraph{}
	Since \(h(\mathcal{K}_1^i) \subseteq \mathcal{K}_2^b\) for each \(b \in \bigcap h(\mathcal{K}_1^i)\),  it suffices to show that for exactly one \(b_i \in \bigcap h(\mathcal{K}_1^i)\) we have:
	\[
		|h(\mathcal{K}_1^i)|=|\mathcal{K}_2^{b_i}|.
	\]
	\paragraph{}
	Assume, for the sake of contradiction, that for all \(b \in \bigcap h(\mathcal{K}_1^i)\),
	we have \(|h(\mathcal{K}_1^i)| < |\mathcal{K}_2^b|\). This implies \(i \notin \bigcap h^{-1}(\mathcal{K}_2^b)\), and thus:
	\[
		\bigcap h^{-1}(\mathcal{K}_2^b) \subset \bigcap \mathcal{K}_1^i.
	\]
	Consequently:
	\[
		\bigcup_{j=1}^n \bigcap h^{-1}(\mathcal{K}_2^{b_j}) \subset \bigcap \mathcal{K}_1^i,
	\]
	which implies:
	\[
		\left| \bigcup_{j=1}^n \bigcap h^{-1}(\mathcal{K}_2^{b_j}) \right| < n.
	\]
	Then by Corollary~\ref{cor:cj}, we have:
	\[
		\left| \bigcup_{j=1}^n \bigcap \mathcal{K}_2^{b_j} \right| < n.
	\]
	However, Corollary~\ref{cor:ui} ensures that:
	\[
		\left| \bigcup_{j=1}^n \bigcap \mathcal{K}_2^{b_j} \right| \geq n,
	\]
	which is a contradiction.
	Therefore, by Lemma~\ref{lem:mi}, there exists exactly one \(j \in \bigcap h(\mathcal{K}_1^i)\) such that:
	\[
		\mathcal{K}_2^j = h(\mathcal{K}_1^i).
	\]
	This completes the proof.
\end{proof}

\begin{theorem}\label{thm:main}
	Let \( \mathcal{K}_1 \) and \( \mathcal{K}_2 \) be two pure union-closed families of sets. For every isomorphism \( h: \mathcal{K}_1 \to \mathcal{K}_2 \), there exists a \textbf{hyperisomorphism} \( H: \bigcup \mathcal{K}_1 \to \bigcup \mathcal{K}_2 \) such that:
	\[
		h(A) = \{ H(a) \mid a \in A \} \quad \text{for all } A \in \mathcal{K}_1.
	\]
\end{theorem}

\begin{proof}
	Let \(i \in \bigcup \mathcal{K}_1\). By Lemma~\ref{lem:mp}, there exists a corresponding \(j \in \bigcup \mathcal{K}_2\) which allows us to define \(H\).
	To show that \(H\) is bijective, suppose \(H(a) = H(b)\) for some \(a, b \in \bigcup \mathcal{K}_1\). This implies that:
	\[
		\mathcal{K}_1^a = \mathcal{K}_1^b.
	\]
	Using Lemma~\ref{lem:mi}, we conclude that \(a = b\), establishing that \(H\) is injective. Furthermore, by Corollary~\ref{cor:ci}, \(H\) is clearly surjective. Thus, \(H\) is bijective.
	Next, let \(A \in \mathcal{K}_1\) and \(i \in A\). To complete the proof, we need to show that \(H(i) \in h(A)\). By Lemma~\ref{lem:mp}, we have
	\[
		h(A) \in \mathcal{K}_2^{H(i)},
	\]
	which satisfies the required condition.
\end{proof}

\newpage

\bibliographystyle{plainurl}
\bibliography{references}

\begin{thebibliography}{10}

\bibitem{AbdollahiWoodroofeZaimi2017}
A.~Abdollahi, R.~Woodroofe, and G.~Zaimi.
\newblock Frankl’s conjecture for subgroup lattices.
\newblock {\em Electron. J. Combin.}, 24(3):Paper P3.25, 9 pp., 2017.
\newblock \href {https://doi.org/10.37236/6248} {\path{doi:10.37236/6248}}.

\bibitem{AbeNakano1998}
T.~Abe and B.~Nakano.
\newblock Frankl’s conjecture is true for modular lattices.
\newblock {\em Graphs Combin.}, 14(4):305--311, 1998.
\newblock \href {https://doi.org/10.1007/PL00021180}
  {\path{doi:10.1007/PL00021180}}.

\bibitem{AbeNakano2000}
T.~Abe and B.~Nakano.
\newblock Lower semimodular types of lattices: {Frankl's} conjecture holds for
  lower quasi-semimodular lattices.
\newblock {\em Graphs Combin.}, 16(1):1--16, 2000.
\newblock \href {https://doi.org/10.1007/s00373-999-0128-5}
  {\path{doi:10.1007/s00373-999-0128-5}}.

\bibitem{Alweiss2024}
R.~Alweiss, B.~Huang, and M.~Sellke.
\newblock Improved lower bound for frankl’s union-closed sets conjecture.
\newblock {\em Electron. J. Combin.}, 31(3):P3.35, 2024.
\newblock \href {https://doi.org/10.37236/12232} {\path{doi:10.37236/12232}}.

\bibitem{Bruhn2015}
H.~Bruhn and O.~Schaudt.
\newblock The journey of the union-closed sets conjecture, 2015.
\newblock URL: \url{https://arxiv.org/abs/1309.3297}, \href
  {http://arxiv.org/abs/1309.3297} {\path{arXiv:1309.3297}}.

\bibitem{Cambie2022}
S.~Cambie.
\newblock Better bounds for the union-closed sets conjecture using the entropy
  approach, 2022.
\newblock URL: \url{https://arxiv.org/abs/2212.12500}, \href
  {http://arxiv.org/abs/2212.12500} {\path{arXiv:2212.12500}}.

\bibitem{CarvalhoMachiavelo2024}
A.~Carvalho and A.~Machiavelo.
\newblock On supratopologies, normalized families and frankl conjecture, 2024.
\newblock URL: \url{https://arxiv.org/abs/2408.11213}, \href
  {http://arxiv.org/abs/2408.11213} {\path{arXiv:2408.11213}}.

\bibitem{Colbert2024}
C.~H. Colbert.
\newblock Chain conditions and optimal elements in generalized union-closed
  families of sets, 2024.
\newblock URL: \url{https://arxiv.org/abs/2412.18740}, \href
  {http://arxiv.org/abs/2412.18740} {\path{arXiv:2412.18740}}.

\bibitem{CzedliSchmidt2008}
G.~Cz{\'e}dli and E.~T. Schmidt.
\newblock Frankl’s conjecture for large semimodular and planar semimodular
  lattices.
\newblock {\em Acta Univ. Palacki. Olomuc., Fac. Rerum Nat. Mathematica},
  47(1):47--53, 2008.
\newblock URL: \url{http://eudml.org/doc/32473}.

\bibitem{DasWu2024}
S.~Das and S.~Wu.
\newblock Frequent elements in union-closed set families, 2024.
\newblock URL: \url{https://arxiv.org/abs/2412.03862}, \href
  {http://arxiv.org/abs/2412.03862} {\path{arXiv:2412.03862}}.

\bibitem{frankl1995}
P.~Frankl.
\newblock Extremal set systems.
\newblock In {\em Handbook of Combinatorics}, volume~2, pages 1293--1329. 1995.

\bibitem{Gendler2025}
G.~Gendler.
\newblock Partial results for union-closed conjectures on the weighted cube,
  2025.
\newblock URL: \url{https://arxiv.org/abs/2504.13347}, \href
  {http://arxiv.org/abs/2504.13347} {\path{arXiv:2504.13347}}.

\bibitem{gilmer2022}
J.~Gilmer.
\newblock A constant lower bound for the union-closed sets conjecture, 2022.
\newblock URL: \url{https://arxiv.org/abs/2211.09055}, \href
  {http://arxiv.org/abs/2211.09055} {\path{arXiv:2211.09055}}.

\bibitem{JoshiWaphare2019}
V.~Joshi and B.~N. Waphare.
\newblock Frankl’s conjecture for breadth two lattices.
\newblock {\em Commun. Algebra}, 47(9):3730--3735, 2019.
\newblock \href {https://doi.org/10.1080/00927872.2019.1593360}
  {\path{doi:10.1080/00927872.2019.1593360}}.

\bibitem{JoshiWaphareKavishwar2016}
V.~Joshi, B.~N. Waphare, and S.~P. Kavishwar.
\newblock A proof of frankl’s union-closed sets conjecture for dismantlable
  lattices.
\newblock {\em Algebra Universalis}, 76(3):351--354, 2016.
\newblock \href {https://doi.org/10.1007/s00012-016-0405-0}
  {\path{doi:10.1007/s00012-016-0405-0}}.

\bibitem{Karpas2017}
I.~Karpas.
\newblock Two results on union-closed families, 2017.
\newblock URL: \url{https://arxiv.org/abs/1708.01434}, \href
  {http://arxiv.org/abs/1708.01434} {\path{arXiv:1708.01434}}.

\bibitem{LuRaz2024}
K.~Lu and A.~Raz.
\newblock A note on the union-closed sets conjecture and reimer's average set
  size theorem, 2024.
\newblock URL: \url{https://arxiv.org/abs/2405.10639}, \href
  {http://arxiv.org/abs/2405.10639} {\path{arXiv:2405.10639}}.

\bibitem{moghaddasmehr2023}
M.~J.~Moghaddas Mehr.
\newblock A note on the union-closed sets conjecture.
\newblock arXiv, 2023.
\newblock \href {http://arxiv.org/abs/2309.01704} {\path{arXiv:2309.01704}},
  \href {https://doi.org/10.48550/arXiv.2309.01704}
  {\path{doi:10.48550/arXiv.2309.01704}}.

\bibitem{Poonen1992}
B.~Poonen.
\newblock Union-closed families.
\newblock {\em J. Combin. Theory Ser. A}, 59(2):253--268, 1992.

\bibitem{Reinhold2000}
J.~Reinhold.
\newblock Frankl’s conjecture is true for lower semimodular lattices.
\newblock {\em Graphs Combin.}, 16(1):115--116, 2000.
\newblock \href {https://doi.org/10.1007/s003730050008}
  {\path{doi:10.1007/s003730050008}}.

\bibitem{Rival1985}
I.~Rival, editor.
\newblock {\em Graphs and Order}, volume 147 of {\em NATO ASI Ser.}
\newblock Springer, Dordrecht, 1985.

\bibitem{Yu2023}
L.~Yu.
\newblock Dimension-free bounds for the union-closed sets conjecture.
\newblock {\em Entropy}, 25(5):767, 2023.
\newblock \href {https://doi.org/10.3390/e25050767}
  {\path{doi:10.3390/e25050767}}.

\end{thebibliography}

\end{document}